\documentclass{article}
\usepackage[round]{natbib}
\usepackage{mathrsfs}
\usepackage{latexsym}
\usepackage{amsmath}
\usepackage{amsfonts}
\usepackage{url}
\usepackage{amssymb}
\usepackage{amsthm}

\newcommand{\inv}{^{-1}} %
\newcommand{\TV}[1]{\nrm{#1}_{\textrm{{\tiny \textup{TV}}}}}
\newcommand{\mix}{_{\textrm{{\tiny \textup{mix}}}}}

\newcommand{\X}{\calX}

\newcommand{\f}{\varphi}
\usepackage{bbm}
\newcommand{\chr}{\boldsymbol{\mathbbm{1}}} %
\newcommand{\pred}[1]{\chr_{\left\{ #1 \right\}}}

\newcommand{\E}{\mathbb{E}}

\newcommand{\lam}{\lambda}

\newcommand{\diam}{\operatorname{diam}}

\renewcommand{\Pr}[1]{\P\!\paren{#1}}
\renewcommand{\P}{\mathbb{P}}

\newcommand{\ssum}[3]{\narsum{ {#1}_{#2}^{#3} }}
\newcommand{\ben}{\begin{enumerate}}
\newcommand{\een}{\end{enumerate}}
\newcommand{\bit}{\begin{itemize}}
\newcommand{\eit}{\end{itemize}}

\def\clap#1{\hbox to 0pt{\hss#1\hss}}

\def\mathclap{\mathpalette\mathclapinternal}

\def\mathclapinternal#1#2{\clap{$\mathsurround=0pt#1{#2}$}}
\newcommand{\narsum}[1]{\sum_{\mathclap{#1}}}

\newcommand{\prd}{^{n}}
\newcommand{\seq}[3]{(#1_{#2},\ldots,#1_{#3})}
\newcommand{\sseq}[3]{#1_{#2}^{#3}}  %

\newcommand{\nrm}[1]{\left\Vert #1 \right\Vert}

\newcommand{\calA}{\mathcal{A}}

\newcommand{\calX}{\mathcal{X}}
\newcommand{\calY}{\mathcal{Y}}
\newcommand{\calZ}{\mathcal{Z}}

\newcommand{\R}{\mathbb{R}}
\newcommand{\N}{\mathbb{N}}

\newcommand{\beq}{\begin{eqnarray*}}
\newcommand{\eeq}{\end{eqnarray*}}
\newcommand{\beqn}{\begin{eqnarray}}
\newcommand{\eeqn}{\end{eqnarray}}
\newcommand{\paren}[1]{\left( #1 \right)}
\newcommand{\sqprn}[1]{\left[ #1 \right]}
\newcommand{\tlprn}[1]{\left\{ #1 \right\}}
\newcommand{\set}[1]{\tlprn{#1}}
\newcommand{\abs}[1]{\left| #1 \right|}

\newcommand{\gn}{\, | \,}

\newcommand{\ts}{\textstyle}
\newcommand{\bs}{\boldsymbol}
\renewcommand{\th}{\ensuremath{^{\mathrm{th}}}~}
\newcommand{\hide}[1]{}
\newcommand{\oo}[1]{\frac{1}{#1}}

\def\eps{\varepsilon}
\newcommand{\defeq}{=}

\newtheorem{theorem}{Theorem}

\newtheorem{lemma}{Lemma}

\newcommand{\bepf}{\begin{proof}}
\newcommand{\enpf}{\end{proof}}

\newcommand{\shortlong}[2]{{#2}}

\title{
Concentration in unbounded metric spaces and algorithmic stability
}

\author{
Aryeh Kontorovich
}

\begin{document}
\maketitle

\begin{abstract}
We prove an extension of McDiarmid's inequality for metric spaces with unbounded diameter.
To this end, 
we introduce the notion of the {\em subgaussian diameter},
which is a 
distribution-dependent
refinement of the metric diameter.
Our technique provides an alternative approach to that of Kutin and Niyogi's 
method of weakly difference-bounded
functions, and yields
nontrivial, dimension-free results in some interesting cases where the former does not.
As an application, we give 
apparently the first
generalization bound in 
the
algorithmic stability setting
that holds for unbounded loss functions.
We give two extensions of the basic concentration result:
to strongly mixing processes
and to other Orlicz norms.
\end{abstract}
\newcommand{\subg}{\Delta_{\textrm{{\tiny \textup{SG}}}}}
\newcommand{\bsubg}{\bar\Delta_{\textrm{{\tiny \textup{SG}}}}}
\newcommand{\orlp}{\Delta_{\textrm{{\tiny \textup{OR}}}(p)}}

\section{Introduction}
Concentration of measure inequalities are at the heart of 
statistical learning theory. Roughly speaking, concentration
allows one to conclude that the performance of a (sufficiently ``stable'')
algorithm on a (sufficiently ``close to iid'') sample is indicative of
the algorithm's performance on future data.
Quantifying 
what it means for an algorithm to be
{\em stable} and for the sampling process to be {\em close to iid}
is by no means straightforward and 
much recent work has been motivated by these questions.
It turns out that the various notions of stability are naturally expressed in terms of
the 
Lipschitz continuity of the algorithm in question
\citep{DBLP:journals/jmlr/BousquetE02a,DBLP:conf/uai/KutinN02,MR2181255,MR2738779}, while 
appropriate relaxations of the 
iid assumption 
are achieved using various kinds of strong mixing
\citep{MR1921877,
gamarnik03,
rostamizadeh07,mohri-rosta08,
DBLP:conf/nips/SteinwartC09,Steinwart2009175,
springerlink:10.1007/s10255-011-0096-4,
DBLP:journals/jmlr/MohriR10,
london:nips12asalsn,london:icml13,
me-cosma-nips}.

An elegant and powerful work-horse driving many of the aforementioned results
is McDiarmid's inequality \citep{mcdiarmid89}:
\beqn
\label{eq:mcd}
\Pr{ 
\abs{
\f
- 
\E \f
}
> t} \le 
2\exp\paren{-\frac{2t^2}{
\sum_{i=1}^nw_i^2
}},
\eeqn
where 
$\f$ is a real-valued function of the sequence of independent random variables $X=(X_1,\ldots,X_n)$, such that
\beqn
\label{eq:lip}
\abs{\f(x)-\f(x')}\le w_i
\eeqn
whenever $x$ and $x'$ differ only in the $i$\th coordinate.
\hide{
$\f$ is $\beta$-Lipschitz with respect to the Hamming metric weighed by $w\in(0,\infty)^n$.
Nor is it possible to obtain concentration bounds that are both distribution-free and do not depend on the
Lipschitz constant of the random variable in question: if each $X_i$ takes the value $\pm R$ independently with
probability $1/2$, then $\sum X_i$ will have typical deviations of $\sqrt n R$.}
Aside from being instrumental in proving PAC bounds \citep{CambridgeJournals:8212764},
McDiarmid's inequality
has also found use in algorithmic stability results \citep{DBLP:journals/jmlr/BousquetE02a}.
Non-iid extensions of (\ref{eq:mcd}) have also been considered \citep{marton96,MR1771956,chazottes07,kontram06}.

The distribution-free nature of 
McDiarmid's inequality
makes it an attractive tool in learning theory, but also
imposes inherent limitations on its applicability. {Chief among these limitations
is the inability of (\ref{eq:mcd}) to provide
risk bounds for unbounded loss functions.}
Even in the bounded case,
if the Lipschitz condition (\ref{eq:lip}) holds
not everywhere but only with high probability --- say, with a much larger constant on a small set of exceptions ---
the bound in (\ref{eq:mcd}) 
still charges the full cost of the worst-case constant. 
To counter this difficulty, \citet{kutin02,DBLP:conf/uai/KutinN02}
introduced an extension of 
McDiarmid's inequality
to {\em weakly difference-bounded} functions and used it to analyze the risk of ``almost-everywhere'' stable
algorithms. This influential result has been invoked in a number of recent papers \citep{MR2277917,Mukherjee06,hush07,Agarwal:2009:GBR:1577069.1577085,MR2738779,6172585}.

However, the approach of \citeauthor{DBLP:conf/uai/KutinN02} entails some difficulties as well.
These come in two flavors: analytical (complex statement and proof) and practical (conditions are still too restrictive in some cases); we will elaborate upon this in Section~\ref{sec:rel-work}. In this paper,
we propose an alternative approach to the concentration of ``almost-everywhere'' or ``average-case'' Lipschitz functions. 
To this end, we introduce the notion of the {\em subgaussian diameter} of a metric probability space.
The latter may be finite even when the metric diameter is infinite, and we show that this notion
generalizes the more restrictive property of bounded differences.

\paragraph{Main results.}
This paper's principal contributions include
defining the subgaussian diameter of a metric probability space
and identifying its role in relaxing the bounded differences condition.
In 
Theorem~\ref{thm:subg-indep},
we show that the subgaussian diameter can essentially replace the far more restrictive metric diameter
in concentration bounds.
This result has direct ramifications for algorithmic stability (Theorem~\ref{thm:alg-stab}).
We furthermore extend our concentration inequality to non-independent processes (Theorem~\ref{thm:subg-mix})
and to other Orlicz norms (Theorem~\ref{thm:p-orl}).

\paragraph{Outline of paper.}
In Section~\ref{sec:def-not} 
we define the subgaussian diameter and relate it to (weakly) bounded differences
in Section~\ref{sec:rel-work}.
We state and prove the concentration inequality based on this notion in Section~\ref{sec:subg-indep}
and give an application to algorithmic stability in Section~\ref{sec:alg-stab}.
We then give an extension to non-independent data in Section~\ref{sec:mix} 
and discuss other Orlicz norms in Section~\ref{sec:other-orlicz}.
Conclusions and some open problems are presented in Section~\ref{sec:open-prob}.

\section{Preliminaries}
\label{sec:def-not} 
A {\em metric probability space} $(\X,\rho,\mu)$ is a measurable space $\X$ whose Borel $\sigma$-algebra is induced by
the metric $\rho$, endowed with the probability measure $\mu$. Our results are most cleanly presented when $\X$
is a discrete set but they continue to hold verbatim for 
Borel probability measures on Polish spaces.
It will be convenient to 
write $\E\f=\sum_{x\in\X}\P(x)\f(x)$ even when the latter is an integral.
Random variables are capitalized ($X$), specified sequences
are written in lowercase,
the 
notation
$\sseq{X}{i}{j}\defeq
\seq{X}{i}{j}
$ is used for all sequences, and
sequence concatenation is denoted multiplicatively:
$\sseq{x}{i}{j}\sseq{x}{j+1}{k}=\sseq{x}{i}{k}$.
We will frequently use the shorthand
$\P(x_i^j)=\prod_{k=i}^j \Pr{X_k=x_k}$.
Standard order of magnitude notation such as $O(\cdot)$ and $\Omega(\cdot)$
will be used.

A function
$\f:\X\to\R$
is $L$-Lipschitz if
\beq
\abs{\f(x)-\f(x')}\le L\rho(x,x'),
\qquad x,x'\in\X.
\eeq

Let $(\X_i,\rho_i,\mu_i)$, $i=1,\ldots,n$  
be a 
sequence of
metric probability spaces. 
We define the product 
probability space
\beq
\X\prd = \X_1\times\X_2\ldots\times\X_n
\eeq
with the product measure\hide{
We will denote partial products by
$$\X_i^j
=\X_i\times\X_{i+1}\times\ldots\times\X_j.$$
}
\beq
\mu\prd=\mu_1\times\mu_2\times\ldots\times\mu_n
\eeq
and $\ell_1$ product
metric 
\beqn
\label{eq:prod-metr}
\rho\prd(x,y) = \sum_{i=1}^n \rho_i(x_i,y_i),\qquad x,y\in\X\prd
.
\eeqn
We will denote partial products by
$$\X_i^j
=\X_i\times\X_{i+1}\times\ldots\times\X_j.$$

We write $X_i\sim\mu_i$ to mean that $X_i$ is an $\X_i$-valued random variable with law $\mu_i$ ---
i.e., $\Pr{X_i\in A}=\mu_i(A)$ for all Borel $A\subset\X_i$.
This notation extends naturally to sequences: $X_1^n\sim\mu\prd$.
We will associate to each 
$(\X_i,\rho_i,\mu_i$)
the {\em symmetrized distance} random variable $\Xi(\X_i)$ defined by
\beqn
\label{eq:zdef}
\Xi(\X_i) = \epsilon_i \rho_i(X_i,x_i'), 
\eeqn
where $X_i,x_i'\sim\mu_i$ are independent
and $\epsilon_i=\pm1$ with probability $1/2$,
independent of $X_i,x_i'$. We note right away that $\Xi(\X_i)$ is a {\em centered}
random variable:
\beqn
\label{eq:zcen}
\E[\Xi(\X_i)]=0.
\eeqn

A real-valued random variable $X$ 
is said to be {\em subgaussian} if it
admits a $\sigma>0$ such that
\beqn
\label{eq:lap}
\E e^{\lam X}\le 
e^{\sigma^2\lam^2/2},
\qquad \lam\in\R.
\eeqn
The smallest $\sigma$ for which (\ref{eq:lap})
holds will be denoted by $\sigma^*(X)$. 
\hide{This characterization has an
(apparently folklore) quantitative version
\citep{Rivasplata2012}:
\begin{lemma}
\label{lem:subg}
Let
$X$ be
a centered subgaussian random variable.
Then 
\begin{itemize}
\item[(i)] 
$\sigma^*(X) \le \sqrt3\nrm{X}_{\psi_2}
;
$
\item[(ii)] 
$
\nrm{X}_{\psi_2}\le
\sqrt 6\sigma^*(X)
$.
\end{itemize}
\end{lemma}
}

We define the 
{\em subgaussian diameter}
$\subg(\X_i)$
of the 
metric probability space
$(\X_i,\rho_i,\mu_i)$ 
in terms 
of its
symmetrized distance
$\Xi(\X_i)$:
\beqn
\label{eq:orldef}
\subg(\X_i) = 
\sigma^*(\Xi(\X_i)).
\eeqn
\hide{
It follows from Lemma~\ref{lem:subg} that
\beq
\nrm{\Xi(\X_i)}_{\psi_2}/\sqrt6
\le
\subg(\X_i) \le \sqrt3\nrm{\Xi(\X_i)}_{\psi_2}.
\eeq
}
If a metric probability space $(\X,\rho,\mu)$ has
finite diameter,
\beq
\diam(\X):=\sup_{x,x'\in\X}\rho(x,x') <\infty,
\eeq
then its subgaussian diameter is also finite:
\begin{lemma}
\label{lem:findim}
\beq
\subg(\X)\le\diam(\X).
\eeq
\end{lemma}
\bepf
Let $\Xi=\Xi(\X)$ be the 
symmetrized distance.
By (\ref{eq:zcen}), we have $\E[\Xi]=0$ and certainly $|\Xi|\le\diam(\X)$. Hence,
\beq
\E e^{\lam \Xi} \le \exp((2\diam(\X)\lam)^2/8)
=
\exp(\diam(\X)^2\lam^2/2),
\eeq
where the inequality follows from Hoeffding's Lemma.
\enpf
The bound in Lemma~\ref{lem:findim} is nearly tight in the sense that for every $\eps>0$ there
is a metric probability space $(\X,\rho,\mu)$ for which
\beqn
\diam(\X)<\subg(\X)+\eps.
\eeqn
To see this, take $\X$ to be an $N$-point space with the uniform distribution
and $\rho(x,x')=1$ for all distinct $x,x'\in\X$. Taking $N$ sufficiently large
makes $\subg(\X)$ arbitrarily close to $\diam(\X)=1$. 
We do not know whether $\diam(\X)=\subg(\X)$ can be achieved.

On the other hand, there exist unbounded metric probability spaces with finite
subgaussian diameter. 
A simple example is 
$(\X,\rho,\mu)$ with
$\X=\R$,
$\rho(x,x')=\abs{x-x'}$
and $\mu$
the standard Gaussian
probability measure $d\mu = (2\pi)^{-1/2}e^{-x^2/2}dx$.
Obviously, $\diam(\X)=\infty$. Now 
the 
symmetrized distance
$\Xi=\Xi(\X)$ 
is distributed as the difference (=sum) of two standard Gaussians:
$\Xi\sim N(0,2)$. Since
$\E e^{\lam \Xi} = e^{\lam^2}$, we have
\beqn
\label{eq:gauss-diam}
\subg(\X)=\sqrt2
.
\eeqn
More generally, the subgaussian
distributions on $\R$ are precisely those for which $\subg(\R)<\infty$.

\section{Related work}
\label{sec:rel-work}
McDiarmid's inequality (\ref{eq:mcd})\hide{
may be stated in the notation of
Section~\ref{sec:def-not} as follows. 
Let $(\X_i,\rho_i,\mu_i)$, $i=1,\ldots,n$  
be a sequence of
metric probability spaces and $\f:\X\prd\to\R$
a $1$-Lipschitz function.
Then
\beqn
\label{eq:mcd2}
\Pr{ 
\abs{
\f
- 
\E \f
}
> t} \le 
2\exp\paren{-\frac{2t^2}{
\sum_{i=1}^n\diam(\X_i)^2
}}.
\eeqn
Being a restatement of (\ref{eq:mcd}),
the bound in (\ref{eq:mcd2})
}
suffers from the limitations mentioned above:
it completely ignores the distribution and is vacuous if even one of
the $w_i$ is infinite.\footnote{Note, though, that McDiarmid's inequality
is sharp in the sense that the constants in (\ref{eq:mcd}) cannot be
improved in a distribution-free fashion.}
In order to address some of these issues, \citet{kutin02,DBLP:conf/uai/KutinN02} proposed 
an
extension of McDiarmid's inequality
to ``almost everywhere'' Lipschitz
functions $\f:\X\prd\to\R$.
To formalize this,
fix an $i\in[n]$ and
let $X_1^n\sim\mu\prd$ and $x_i'\sim\mu_i$ be independent.
Define $\tilde X_1^n=\tilde X_1^n(i)$ by
\beqn
\label{eq:xtil}
\tilde X_j(i) = \left\{
\begin{array}{ll}
X_j, & j\neq i \\
x_i', & j=i.
\end{array}
\right.
\eeqn
\citeauthor{DBLP:conf/uai/KutinN02} define $\f$ to be 
{\em weakly difference-bounded} 
by $(b,c,\delta)$ if
\beqn
\label{eq:k-n-b}
\Pr{|\f(X)-\f(\tilde X(i))|>b}=0
\eeqn
and
\beqn
\label{eq:k-n-c}
\Pr{|\f(X)-\f(\tilde X(i))|>c}<\delta
\eeqn
for all $1\le i\le n$.

The precise 
result 
of
\citet[Theorem 1.10]{kutin02}
is somewhat unwieldy
to state --- indeed, the present work was motivated in part
by a desire for simpler tools.
Assuming that 
$\f$ is weakly difference-bounded
by $(b,c,\delta)$ 
with
\beqn
\label{eq:delta-ass}
\delta=\exp(-\Omega(n))
\eeqn
and $c=O(1/n)$,
their bound states that
\beqn
\label{eq:k-n-ineq}
\Pr{\abs{\f-\E\f}\ge t}\le \exp(-\Omega(nt^2))
\eeqn
for a certain range of $t$ and $n$.
As noted by \citet{MR2181255}, the exponential decay assumption (\ref{eq:delta-ass})
is necessary in order for the Kutin-Niyogi method to yield exponential concentration.
In contrast, the bounds we prove here 
\begin{itemize}
\item[(i)] do not require $|\f(X)-\f(\tilde X)|$ to be everywhere bounded as in (\ref{eq:k-n-b})
\item[(ii)] have a simple statement and proof, and generalize to non-iid processes with relative ease.
\end{itemize}
We defer the quantitative comparisons between (\ref{eq:k-n-ineq}) and our results until the latter are formally stated
in Section~\ref{sec:subg-indep}.

In a different line of work, \citet{MR2425108} considered an extension of Hoeffding's inequality
to unbounded random variables. His bound only holds for sums (as opposed to general Lipschitz functions)
and the summands must be non-negative (i.e., unbounded only in the positive direction).
An earlier notion of ``effective'' metric diameter in the context of concentration is that of {\em metric space length} \citep{schechtman82}.
Another distribution-dependent refinement of diameter is the {\em spread
constant} \citep{along-boppana-spencer-98}.
\citet{lecue2013} gave minimax bounds for empirical risk minimization over subgaussian classes.

\section{Concentration via subgaussian diameter}
\label{sec:subg-indep}
McDiarmid's inequality (\ref{eq:mcd}) 
may be stated in the notation of
Section~\ref{sec:def-not} as follows. 
Let $(\X_i,\rho_i,\mu_i)$, $i=1,\ldots,n$  
be a sequence of
metric probability spaces and $\f:\X\prd\to\R$
a $1$-Lipschitz function.
Then
\beqn
\label{eq:mcd2}
\Pr{ 
\abs{
\f
- 
\E \f
}
> t} \le 
2\exp\paren{-\frac{2t^2}{
\sum_{i=1}^n\diam(\X_i)^2
}}.
\eeqn
We defined the subgaussian diameter $\subg(\X_i)$
in Section~\ref{sec:def-not}, showing in Lemma~\ref{lem:findim}
that it never exceeds the metric diameter. We also showed by
example that the former can be finite when the latter is infinite.
The main result of this section is that $\diam(\X_i)$ 
in (\ref{eq:mcd2})
can essentially
be replaced by $\subg(\X_i$):
\begin{theorem}
\label{thm:subg-indep}
If 
$\f:\X\prd\to\R$ is $1$-Lipschitz then
$\E\f<\infty$ and
\beq
\Pr{ 
\abs{
\f
- 
\E \f
}
> t} \le 
2\exp\paren{-\frac{t^2}{
2
\sum_{i=1}^n\subg^2(\X_i)
}}.
\eeq
\end{theorem}
Our constant in the exponent is worse than that of (\ref{eq:mcd2}) by a factor of $4$;
this appears to be an inherent artifact of our method.
\bepf
The strong integrability of $\f$ --- and in particular, finiteness of $\E\f$ --- follow from
exponential concentration \citep{ledoux01}. The rest of the proof will proceed via the 
Azuma-Hoeffding-McDiarmid method of martingale differences.
Define $V_i=\E[\f\gn X_1^i] - \E[\f\gn X_1^{i-1}]$ and expand
\beq
\E[\f\gn X_1^i] &=& \sum_{x_{i+1}^n\in\X_{i+1}^n}\P(x_{i+1}^n)\f(X_1^i x_{i+1}^n)\\
\E[\f\gn X_1^{i-1}] &=& \sum_{x_{i}^n\in\X_{i}^n}\P(x_{i}^n)\f(X_1^{i-1} x_{i}^n).
\eeq
Let $\tilde V_i$ be $V_i$ conditioned on $\sseq X1{i-1}$;
thus,
\shortlong{
\begin{align*}
&\tilde V_i=
\sum_{x_{i+1}^n}
\P(x_{i+1}^n)
\cdot\\&
\sum_{x_i,x_i'
}
\P(x_i)\P(x_i')
\paren{
\f(\sseq X1{i-1}x_i\sseq{x}{i+1}{n}) - \f(\sseq X1{i-1}x_i'\sseq x{i+1}n)
}.
\end{align*}
}{
\begin{align*}
&\tilde V_i=
\sum_{x_{i+1}^n}
\P(x_{i+1}^n)
\sum_{x_i,x_i'
}
\P(x_i)\P(x_i')
\paren{
\f(\sseq X1{i-1}x_i\sseq{x}{i+1}{n}) - \f(\sseq X1{i-1}x_i'\sseq x{i+1}n)
}.
\end{align*}
}

Hence, by Jensen's inequality, we have
\shortlong{
\begin{align*}
&\E[e^{ \lam V_i}\gn X_1^{i-1}] \le
\sum_{x_{i+1}^n
}
\P(\sseq x{i+1}n)
\cdot\\
&
\narsum{y,y'
} \P(y)\P(y') 
e^{\lam(
\f(X_1^{i-1}yx_{i+1}^n)
-
\f(X_1^{i-1}y'x_{i+1}^n)
)}.
\end{align*}
}{
\begin{align*}
&\E[e^{ \lam V_i}\gn X_1^{i-1}] \le
\sum_{x_{i+1}^n
}
\P(\sseq x{i+1}n)
\narsum{y,y'
} \P(y)\P(y') 
e^{\lam(
\f(X_1^{i-1}yx_{i+1}^n)
-
\f(X_1^{i-1}y'x_{i+1}^n)
)}.
\end{align*}
}
For fixed 
$
X_1^{i-1}
\in \sseq\X1{i-1}$
and
$
x_{i+1}^{n}
\in \sseq\X{i+1}{n}$,
define $F:\X_i\to\R$ by 
$F(y)=\f(X_1^{i-1}y x_{i+1}^{n})$,
and observe that $F$ is $1$-Lipschitz with respect to $\rho_i$. Since
$e^t+e^{-t}=2\cosh(t)$ and
$\cosh(t)\le\cosh(s)$
for all $|t|\le s$, 
we have\footnote{
An analogous symmetrization technique is employed in \url{http://terrytao.wordpress.com/2009/06/09/talagrands-concentration-inequality}
as a variant of the ``square and rearrange'' trick.
}
\beq
e^{\lam(F(y)-F(y'))}
+
e^{\lam(F(y')-F(y))}
\le
e^{\lam \rho_i(y,y')}+e^{-\lam \rho_i(y,y')},
\eeq
and hence
\shortlong{
\begin{align}
\label{eq:symm-arg}
&\sum_{y,y'\in\X_i} 
\P(y)\P(y') 
e^{\lam(
F(y)
-
F(y')
)}
\\&\le\nonumber
{\ts\oo2}
\!\sqprn{
\narsum{y,y'
} 
\P(y)\P(y') 
e^{\lam \rho_i(y,y')}+
\narsum{
y,y'
} 
\P(y)\P(y') 
e^{-\lam \rho_i(y,y')}
}
\\&=\nonumber
\E e^{\lam \Xi(\X_i)} \le \exp(\lam^2\subg^2(\X_i)/2),
\end{align}
}{
\begin{align}
\label{eq:symm-arg}
&\sum_{y,y'\in\X_i} 
\P(y)\P(y') 
e^{\lam(
F(y)
-
F(y')
)}
\\&\le\nonumber
{\ts\oo2}
\sqprn{
\narsum{y,y'
} 
\P(y)\P(y') 
e^{\lam \rho_i(y,y')}+
\narsum{
y,y'
} 
\P(y)\P(y') 
e^{-\lam \rho_i(y,y')}
}
\\&=\nonumber
\E e^{\lam \Xi(\X_i)} \le \exp(\lam^2\subg^2(\X_i)/2),
\end{align}
}where $\Xi(\X_i)$ is the symmetrized distance (\ref{eq:zdef})
and the last inequality holds by
definition of subgaussian diameter (\ref{eq:lap},\ref{eq:orldef}).
It follows that
\beqn
\E[ e^{\lam V_i} \gn X_1^{i-1}]\le \exp(\lam^2\subg^2(\X_i)/2).
\eeqn
Applying the standard 
Markov's inequality and exponential bounding
argument, we have
\shortlong{
\begin{align}
\nonumber
&\Pr{ \f - \E \f > t} =
\Pr{ \sum_{i=1}^n V_i  > t} 
\\\nonumber&\le
e^{-\lam t}
\E\sqprn{\prod_{i=1}^n e^{\lam V_i}} 
\\\nonumber&= 
e^{-\lam t}
\E\sqprn{\prod_{i=1}^n \E[e^{\lam V_i}\gn X_1^{i-1}]} 
\\\nonumber&\le
e^{-\lam t}
\E\sqprn{\prod_{i=1}^n 
\exp(\lam^2\subg^2(\X_i)/2)
} 
\\&=
\label{eq:exp-bd}
\exp\paren{\oo2\lam^2\sum_{i=1}^n\subg^2(\X_i)-\lam t}.
\end{align}
}{
\begin{align}
\nonumber
\Pr{ \f - \E \f > t} &=
\Pr{ \sum_{i=1}^n V_i  > t} 
\\\nonumber&\le
e^{-\lam t}
\E\sqprn{\prod_{i=1}^n e^{\lam V_i}} 
\\\nonumber&= 
e^{-\lam t}
\E\sqprn{\prod_{i=1}^n \E[e^{\lam V_i}\gn X_1^{i-1}]} 
\\\nonumber&\le
e^{-\lam t}
\E\sqprn{\prod_{i=1}^n 
\exp(\lam^2\subg^2(\X_i)/2)
} 
\\&=
\label{eq:exp-bd}
\exp\paren{\oo2\lam^2\sum_{i=1}^n\subg^2(\X_i)-\lam t}.
\end{align}
}Optimizing over $\lam$ and applying the same argument to $\E\f-\f$ yields our claim.
\enpf

Let us see how
Theorem~\ref{thm:subg-indep} compares to previous results
on some examples.
Consider $\R^n$ equipped with the $\ell_1$ metric $\rho^n(x,x')=\sum_{i\in[n]}|x_i-x_i'|$
and the standard Gaussian product measure $\mu^n=N(0,I_n)$. Let $\f:\R^n\to\R$
be $1/n$-Lipschitz. Then 
Theorem~\ref{thm:subg-indep} yields
(recalling the calculation in (\ref{eq:gauss-diam}))
\beqn
\label{eq:ex1}
\Pr{\abs{\f-\E\f}>\eps}\le2\exp(-n\eps^2/4),\qquad \eps>0,
\eeqn
whereas the inequalities of
McDiarmid (\ref{eq:mcd}) and
Kutin-Niuyogi (\ref{eq:k-n-ineq})
are both uninformative since the metric diameter is infinite.

For our next example, fix an $n\in\N$ and put $\X_i=\set{\pm1,\pm n}$
with the metric $\rho_i(x,x')=|x-x'|$ and the 
distribution $\mu_i(x)\propto e^{-x^2}$.
One may verify via a calculation analogous to (\ref{eq:gauss-diam}) that
$\subg(\X_i)\le\sqrt2$. For independent $X_i\sim\mu_i$, $i=1,\ldots,n$,
put $\f(X_1^n)=n\inv\sum_{i=1}^n X_i$. Then 
Theorem~\ref{thm:subg-indep}
implies that in this case the bound in (\ref{eq:ex1}) holds verbatim.
On the other hand, $\f$ is easily seen to be 
weakly difference-bounded
by $(1,1/n,e^{-\Omega(n)})$ and thus (\ref{eq:k-n-ineq}) also yields
subgaussian concentration, albeit with worse constants. 
Applying (\ref{eq:mcd})
yields the much cruder estimate 
\beq
\Pr{\abs{\f-\E\f}>\eps}\le2\exp(-2\eps^2).
\eeq

\section{Application to algorithmic stability}
\label{sec:alg-stab}
We refer the reader to 
\citep{DBLP:journals/jmlr/BousquetE02a,DBLP:conf/uai/KutinN02,MR2181255} for background on algorithmic
stability and supervised learning. 
Our metric probability space $(\calZ_i,\rho_i,\mu_i)$ will now have the structure $\calZ_i=\calX_i\times\calY_i$
where $\calX_i$ 
and
$\calY_i$ are, respectively, the {\em instance} and {\em label} space of the $i$\th example.
Under the iid assumption, the
$(\calZ_i,\rho_i,\mu_i)$ are identical for all $i\in\N$ 
(and so we will henceforth drop the subscript $i$ from these).
A training sample is 
$S=Z_1^n\sim\mu^n$
is drawn and a {\em learning algorithm} 
$\calA
$ 
inputs $S$ and outputs
a {\em hypothesis} $f:\calX\to\calY$.
The hypothesis 
$f=\calA(S)$
will be denoted by $\calA_S$.
In line with the previous literature, we assume that $\calA$ is symmetric
(i.e., invariant under permutations of $S$).
The {\em loss} of a hypothesis $f$ on an example $z=(x,y)$ is defined 
by 
$$L(f,z)=\ell(f(x),y),$$ where $\ell:\calY\times\calY\to[0,\infty)$ 
is the {\em cost
function}. 
To our knowledge, all previous work required the loss to be bounded 
by some constant $M<\infty$, which figures explicitly in the bounds;
we make no such restriction.

In the algorithmic stability setting, the {\em empirical risk} $\hat R_n(\calA,S)$
is typically defined as
\beqn
\label{eq:emp-risk}
\hat R_n(\calA,S) = \oo n\sum_{i=1}^nL(\calA_S,z_i)
\eeqn
and the {\em true risk} $R(\calA,S)$ as
\beqn
\label{eq:tru-risk}
R(\calA,S) = 
\E_{z\sim\mu}[L(\calA_S,z)].
\eeqn
The goal is to bound the {\em excess risk} $R(\calA,S)-\hat R_n(\calA,S)$.
To this end, a myriad of notions of hypothesis stability have been proposed.
A variant of
{\em uniform stability} in the sense of \citet{MR2181255} ---
which is slightly more general than the homonymous notion
in \citet{DBLP:journals/jmlr/BousquetE02a}
--- may be defined as follows.
The algorithm $\calA$ is said to be $\beta$-uniform stable if
for all $\tilde z\in\calZ$, the function $\f_{\tilde z}:\calZ^n\to\R$ given by
$\f_{\tilde z}(z)=L(\calA_z,\tilde z)$ is $\beta$-Lipschitz with respect to
the Hamming metric on $\calZ^n$:
\shortlong{
\beq
\forall \tilde z\in\calZ,
\forall
z,z'\in\calZ^n:
~
|\f_{\tilde z}(z)-\f_{\tilde z}(z')|\le
\beta\!\sum_{i=1}^n
\!\pred{z_i\neq z'_i}.
\eeq
}{
\beq
\forall \tilde z\in\calZ,
\forall
z,z'\in\calZ^n:
~
|\f_{\tilde z}(z)-\f_{\tilde z}(z')|\le
\beta\sum_{i=1}^n
\pred{z_i\neq z'_i}.
\eeq
}We define the algorithm $\calA$ to be $\beta$-{\em totally Lipschitz stable}
if
the function $\f:\calZ^{n+1}\to\R$ given by
$\f(z_1^{n+1})=L(\calA_{z_1^n},z_{n+1})$ is $\beta$-Lipschitz with respect to
the $\ell_1$ product metric on $\calZ^{n+1}$:
\shortlong{
\beqn
\label{eq:tot-lip}
\forall
z,z'\in\calZ^{n+1}:
~
|\f(z)-\f(z')|\le\beta\!\sum_{i=1}^{n+1}\rho(z_i, z'_i).
\eeqn
}{
\beqn
\label{eq:tot-lip}
\forall
z,z'\in\calZ^{n+1}:
~
|\f(z)-\f(z')|\le\beta\sum_{i=1}^{n+1}\rho(z_i, z'_i).
\eeqn
}Note that total Lipschitz stability is stronger than
uniform stability since it requires the algorithm to
respect the metric of $\calZ$.

Let us bound the bias of stable algorithms.
\begin{lemma}
\label{lem:bias}
Suppose $\calA$ is a
symmetric, 
$\beta$-{totally Lipschitz stable}
learning algorithm
over the metric probability space $(\calZ,\rho,\mu)$
with $\subg(\calZ)<\infty$. Then
\beq
\E[ R(\calA,S)-\hat R_n(\calA,S) ] &\le &
{\ts\oo2}\beta^2\subg^2(\calZ).
\eeq
\end{lemma}
\bepf
Observe, as in the proof of \citep[Lemma 7]{DBLP:journals/jmlr/BousquetE02a},
that for all $i\in[n]$,
\shortlong{
\begin{align}
&\label{eq:bous-decomp}
\E[ R(\calA,S)-\hat R_n(\calA,S) ]
=\\
&\nonumber
\E_{Z_1^n,\tilde Z_1^n}[L(\calA_{Z_1^n},\tilde Z_i) - L(\calA_{\tilde Z_1^n},\tilde Z_i)],
\end{align}
}{
\begin{align}
&\label{eq:bous-decomp}
\E[ R(\calA,S)-\hat R_n(\calA,S) ]
=
\E_{Z_1^n,\tilde Z_1^n}[L(\calA_{Z_1^n},\tilde Z_i) - L(\calA_{\tilde Z_1^n},\tilde Z_i)],
\end{align}
}where $Z_1^n\sim\mu^n$ and $\tilde Z$ is generated from $Z$ via the process
defined in (\ref{eq:xtil}). For fixed $i\in[n]$ and $Z_1^{i-1}$,$Z_{i+1}^n$, define
\beq
W_i(Z_i,Z'_i) = L(\calA_{Z_1^n},Z'_i) - L(\calA_{Z_1^{i-1}Z_i'Z_{i+1}^n},Z'_i)
\eeq
and note that (\ref{eq:tot-lip}) implies 
that $|W_i(Z_i,Z_i')|\le\beta\rho(Z_i,Z_i')$.
Now rewrite (\ref{eq:bous-decomp}) as
\shortlong{
\begin{align}
\label{eq:bias-rewrite}
&\E[ R(\calA,S)-\hat R_n(\calA,S) ] = 
\\&\nonumber
\narsum{z_1^{i-1},z_{i+1}^n}
\P(z_1^{i-1})\P(z_{i+1}^n)
\narsum{z_i,z_i'}\P(z_i)\P(z_i')W_i(z_i,z_i').
\end{align}
}{
\begin{align}
\label{eq:bias-rewrite}
&\E[ R(\calA,S)-\hat R_n(\calA,S) ] = 
\narsum{z_1^{i-1},z_{i+1}^n}
\P(z_1^{i-1})\P(z_{i+1}^n)
\narsum{z_i,z_i'}\P(z_i)\P(z_i')W_i(z_i,z_i').
\end{align}
}Invoking Jensen's inequality and the argument in (\ref{eq:symm-arg}),
\shortlong{
\begin{align*}
&\exp\!\!\paren{\narsum{z_i,z_i'}\P(z_i)\P(z_i')W_i(z_i,z_i')\!}\!
\le \narsum{z_i,z_i'}\P(z_i)\P(z_i')e^{W_i(z_i,z_i')}
\\&=
{\ts\oo2}
\!\sqprn{
\narsum{z_i,z_i'}\P(z_i)\P(z_i')e^{W_i(z_i,z_i')}
+
\narsum{z_i,z_i'}\P(z_i)\P(z_i')e^{-W_i(z_i,z_i')}
}
\\&\le
{\ts\oo2}
\!\sqprn{
\narsum{z_i,z_i'}\P(z_i)\P(z_i')e^{\beta\rho(z_i,z_i')}
+
\narsum{z_i,z_i'}\P(z_i)\P(z_i')e^{-\beta\rho(z_i,z_i')}
}
\\&\le
\exp({\ts\oo2}\beta^2\subg^2(\calZ)).
\end{align*}
}{
\begin{align*}
&\exp\paren{\narsum{z_i,z_i'}\P(z_i)\P(z_i')W_i(z_i,z_i')}
\le \narsum{z_i,z_i'}\P(z_i)\P(z_i')e^{W_i(z_i,z_i')}
\\&=
{\ts\oo2}
\sqprn{
\narsum{z_i,z_i'}\P(z_i)\P(z_i')e^{W_i(z_i,z_i')}
+
\narsum{z_i,z_i'}\P(z_i)\P(z_i')e^{-W_i(z_i,z_i')}
}
\\&\le
{\ts\oo2}
\sqprn{
\narsum{z_i,z_i'}\P(z_i)\P(z_i')e^{\beta\rho(z_i,z_i')}
+
\narsum{z_i,z_i'}\P(z_i)\P(z_i')e^{-\beta\rho(z_i,z_i')}
}
\\&\le
\exp({\ts\oo2}\beta^2\subg^2(\calZ)).
\end{align*}
}
Taking logarithms yields the estimate
\beqn
\label{eq:W-bd}
{\narsum{z_i,z_i'}\P(z_i)\P(z_i')W_i(z_i,z_i')}
\le
{\ts\oo2}\beta^2\subg^2(\calZ),
\eeqn
which, after substituting (\ref{eq:W-bd}) into (\ref{eq:bias-rewrite}),
proves the claim.
\enpf

We now turn to the Lipschitz continuity of the excess risk.
\begin{lemma}
\label{lem:risk-lip}
Suppose $\calA$ is a
symmetric, 
$\beta$-{totally Lipschitz stable}
learning algorithm
and
define the excess risk function
$\f:\calZ^n\to\R$ by $\f(z)=R(\calA,z)-\hat R_n(\calA,z)$.
Then $\f$ is $3\beta$-Lipschitz.
\end{lemma}
\bepf
We examine the two summands separately.
The definition (\ref{eq:tru-risk}) of $R(\calA,\cdot)$
implies that the latter is $\beta$-Lipchitz since
it is a convex combination of $\beta$-Lipschitz functions.
Now $\hat R_n(\calA,\cdot)$ defined in (\ref{eq:emp-risk})
is also a convex combination of $\beta$-Lipschitz functions,
but because $z_i$ appears twice in $L(\calA_{z_1^n},z_i)$,
changing $z_i$ to $z_i'$ could incur a difference of up to
$2\beta\rho(z_i,z_i')$. Hence, $\hat R_n(\calA,\cdot)$ is
$2\beta$-Lipschitz. As Lipschitz constants are subadditive, 
the claim is proved.
\enpf

Combining Lemmas~\ref{lem:bias} and \ref{lem:risk-lip}
with our concentration inequality in Theorem~\ref{thm:subg-indep}
yields the main result of this section:
\begin{theorem}
\label{thm:alg-stab}
Suppose $\calA$ is a
symmetric, 
$\beta$-{totally Lipschitz stable}
learning algorithm
over the metric probability space $(\calZ,\rho,\mu)$
with $\subg(\calZ)<\infty$.
Then, for training samples $S\sim\mu^n$ and $\eps>0$, we have
\shortlong{
\begin{align*}
&\Pr{ 
R(\calA,S)-\hat R_n(\calA,S)>{\ts\oo2}\beta^2\subg^2(\calZ)+\eps
}
\\&\le
\exp\paren{ -\frac{\eps^2}
{
18\beta^2\subg^2(\calZ)n
}}
.
\end{align*}
}{
\begin{align*}
&\Pr{ 
R(\calA,S)-\hat R_n(\calA,S)>{\ts\oo2}\beta^2\subg^2(\calZ)+\eps
}
\le
\exp\paren{ -\frac{\eps^2}
{
18\beta^2\subg^2(\calZ)n
}}
.
\end{align*}
}
\end{theorem}
As in \citet{DBLP:journals/jmlr/BousquetE02a} and related results on algorithmic
stability, we require $\beta=O(1/n)$ for exponential decay.
\citeauthor{DBLP:journals/jmlr/BousquetE02a} showed that this is indeed the case
for some popular learning algorithms, albeit in their less restrictive
definition of stability. We conjecture that many of these algorithms continue
to be stable in our stronger sense and plan to explore this in future work.

\section{Relaxing the independence assumption}
\label{sec:mix}
In this section we generalize Theorem~\ref{thm:subg-indep} to strongly mixing processes.
To this end, we require
some standard facts concerning
the probability-theoretic notions of {\em coupling} 
and
{\em transportation}
\citep{lindvall02,MR1964483,MR2459454}.
Given the probability measures $\mu,\mu'$ on a measurable space $\X$,
a coupling $\pi$ of $\mu,\mu'$ is any probability measure on $\X\times\X$ with
marginals $\mu$ and $\mu'$, respectively. Denoting by $\Pi=\Pi(\mu,\mu')$ the set of
all couplings, we have
\begin{align}
\nonumber
\inf_{\pi\in\Pi}
\pi(\set{(x,y)\in\X^2:x\neq y})&={\ts\oo2}\sum_{x\in\X}|\mu(x)-\mu'(x)|\\
\label{eq:tv}
&=\TV{\mu-\mu'}
\end{align}
where $\TV{\cdot}$ is the {\em total variation} norm. 
An {\em optimal} coupling 
is one that achieves the infimum in (\ref{eq:tv}); one always exists, though it may not be unique.
Another elementary property of couplings is that 
for any two $f,g:\X\to\R$
and any
coupling $\pi\in\Pi(\mu,\mu')$, we have
\beqn
\label{eq:fg}
\E_\mu f-\E_{\mu'}g = \E_{(X,X')\sim\pi}[f(X)-g(X')].
\eeqn

It is possible to refine
the total variation distance between $\mu$ and $\mu'$ so as to respect
the metric of $\X$. Given a space
equipped with
probability measures $\mu,\mu'$ and metric $\rho$, define
the {\em transportation cost\footnote{
This fundamental notion is also known as the 
{\em Wasserstein,
Monge-Kantorovich}, or {\em earthmover} distance;
see \citet{MR1964483,MR2459454} for an encyclopedic treatment.
The use of coupling and transportation techniques to obtain concentration for dependent random variables goes back to
\citet{marton96,samson00,chazottes07}.
} distance} $T_\rho(\mu,\mu')$ by
\beq
T_\rho(\mu,\mu') = \inf_{\pi\in\Pi(\mu,\mu')}\E_{(X,X')\sim\pi} \rho(X,X').
\eeq
It is easy to verify that $T_\rho$ is a valid metric on 
probability measures
and that for $\rho(x,x')=\pred{x\neq x'}$, we have $T_\rho(\mu,\mu')=\TV{\mu-\mu'}$.

As in Section~\ref{sec:subg-indep}, we consider a sequence of metric spaces
$(\X_i,\rho_i)$, $i=1,\ldots,n$ and their $\ell_1$ product $(\X^n,\rho^n)$.
Unlike the independent case, we will allow nonproduct probability measures $\nu$ on $(\X^n,\rho^n)$.
We will write $X_1^n\sim\nu$ to mean that $\Pr{X_1^n\in A}=\nu(A)$ for all Borel $A\subset\X^n$.
For $1\le i\le j<k\le l\le n$, we will use the shorthand
\beq
\P(x_k^l \gn x_i^j) = \Pr{X_k^l=x_k^l \gn X_i^j=x_i^j}.
\eeq
The notation $\P(X_i^j)$ means the marginal distribution of $X_i^j$.
Similarly, $\P(X_k^l \gn X_i^j = x_i^j)$ will denote the conditional distribution.
For $1\le i <n$, and $x_1^i\in\X_1^i$, $x_i'\in\X_i$ define 
\shortlong{
\beq
\tau_{i}(x_1^i,\!x_i') \!=\! T_{\rho_{i+1}^n}
\!(\P(X_{i+1}^n| X_1^i\!=\!x_1^i),\!\P(X_{i+1}^n | X_1^i\!=\!x_1^{i-1}\!x_i')),
\eeq
}{
\beq
\tau_{i}(x_1^i,x_i') = T_{\rho_{i+1}^n}
(\P(X_{i+1}^n\gn X_1^i=x_1^i),\P(X_{i+1}^n \gn X_1^i=x_1^{i-1}x_i')),
\eeq
}where $\rho_{i+1}^n$ is the $\ell_1$ product of $\rho_{i+1},\ldots\rho_n$ as in (\ref{eq:prod-metr}),
and
\beq
\bar\tau_{i} = \sup_{x_1^i\in\X_1^i, x_i'\in\X_i} \tau_{i}(x_1^i,x_i'),
\eeq
with $\bar\tau_n\equiv0$.
In words, $\tau_{i}(x_1^i,x_i')$ measures the transportation cost distance
between the conditional distributions induced on the ``tail'' $\X_{i+1}^n$ 
given two prefixes that differ in the $i$\th coordinate,
and $\bar\tau_{i}$ is the maximal value of this quantity.
\citet{kont07-thesis,kontram06} discuss how to handle
conditioning on measure-zero sets and other technicalities.
Note that for product measures the conditional distributions
are identical and hence $\bar\tau_{i}=0$.

We need one more definition before stating our main result.
For the prefix $x_1^{i-1}$, define the conditional distribution
$$\nu_i(x_1^{i-1})=\Pr{X_i\gn X_1^{i-1}= x_1^{i-1}}$$
and consider the corresponding metric probability space $(\X_i,\rho_i,\nu_i(x_1^{i-1}))$.
Define its {\em conditional subgaussian diameter} by
\beq
\subg(\X_i\gn x_1^{i-1}) = \subg(\X_i,\rho_i,\nu_i(x_1^{i-1}))
\eeq
and the {\em maximal subgaussian diameter} by
\beqn
\label{eq:bsubg}
\bsubg(\X_i) = \sup_{x_1^{i-1}\in\X_1^{i-1}} 
\subg(\X_i\gn x_1^{i-1}).
\eeqn
\hide{
Let us arrange the mixing coefficients $\bar\eta_{ij}$
into the upper-triangular $n\times n$ matrix $H\mix$,
\beq
(H\mix)_{ij} = \left\{
\begin{array}{ll}
\bar\eta_{ij}, & i<j \\
1, & i=j \\
0, & i>j
\end{array}
\right.
\eeq
and
the maximal subgaussian diameters into the vector $\dsubg\in\R^n$,
\beqn
\dsubg=( \bsubg(\X_1),\bsubg(\X_2),\ldots,\bsubg(\X_n)).
\eeqn
}
Note that for product measures, 
(\ref{eq:bsubg}) reduces to
the former definition (\ref{eq:orldef}).
With these definitions, we may state the main result of this section.
\begin{theorem}
\label{thm:subg-mix}
If 
$\f:\X^n\to\R$ is $1$-Lipschitz with respect to $\rho^n$,
then
\beq
\Pr{ 
\abs{
\f
- 
\E \f
}
> t} \le 
2\exp\paren{-\frac{(t-\sum_{i\le n}\bar\tau_i)^2}{
2
\sum_{i\le n} \bsubg^2(\X_i)
}},
\quad t>0.
\eeq
\end{theorem}
Observe that we recover Theorem~\ref{thm:subg-indep} as a special case.
Since typically we will take $t=\eps n$, it suffices that
$\sum_{i\le n}\bar\tau_i=o(n)$ and $\sum_{i\le n}\bsubg^2(\X_i)=O(n)$ 
to ensure a exponential bound
with decay rate $\exp(-\Omega(n\eps^2))$.
\bepf
We begin by considering the martingale difference
$$V_i=\E[\f\gn X_1^i=x_1^i] - \E[\f\gn X_1^{i-1}=x_1^{i-1}]$$
as in the proof of Theorem~\ref{thm:subg-indep}.
More explicitly,
\shortlong{
\begin{align*}
&V_i = 
\\&
\ssum{x}{i+1}{n}\P(x_{i+1}^n| x_1^{i})\f(x_1^ix_{i+1}^n)
-
\ssum{x}{i}{n}\P(x_{i}^n| x_1^{i-1})\f(x_1^{i-1}x_{i}^n)
\\&
=\narsum{x_i'} \P(x_i'\gn x_1^{i-1})
\cdot
\ssum{x}{i+1}{n}
\big[
\P(x_{i+1}^n| x_1^{i})\f(x_1^ix_{i+1}^n)
-
\\&
\phantom{\narsum{x_i'} \P(x_i'\gn x_1^{i-1})\cdot\ssum{x}{i+1}{n}}
\P(x_{i+1}^n| x_1^{i-1} x_i')\f(x_1^{i-1}x_i' x_{i+1}^n)
\big].
\end{align*}
}{
\begin{align*}
V_i &= 
\ssum{x}{i+1}{n}\P(x_{i+1}^n\gn x_1^{i})\f(x_1^ix_{i+1}^n)
-
\ssum{x}{i}{n}\P(x_{i}^n\gn x_1^{i-1})\f(x_1^{i-1}x_{i}^n)
\\&
=\narsum{x_i'} \P(x_i'\gn x_1^{i-1})
\ssum{x}{i+1}{n}
\big[
\P(x_{i+1}^n\gn x_1^{i})\f(x_1^ix_{i+1}^n)
-
\P(x_{i+1}^n\gn x_1^{i-1} x_i')\f(x_1^{i-1}x_i' x_{i+1}^n)
\big].
\end{align*}
}
Define $\tilde V_i$ to be $V_i$ conditioned on $X_1^{i-1}$.
Then
\shortlong{
\begin{align}
\label{eq:Vtil}
&\tilde V_i=
\sum_{x_i,x_i'} 
\P( x_i\gn X_1^{i-1})
\P(x_i'\gn X_1^{i-1})\cdot
\\\nonumber&
\ssum{x}{i+1}{n}
[
\P(x_{i+1}^n\gn X_1^{i-1} x_i)\f(X_1^{i-1}x_i' x_{i+1}^n)
\\\nonumber
&-
\P(x_{i+1}^n\gn X_1^{i-1} x_i')\f(X_1^{i-1}x_i x_{i+1}^n)
].
\end{align}
}{
\begin{align}
\label{eq:Vtil}
&\tilde V_i=
\sum_{x_i,x_i'} 
\P( x_i\gn X_1^{i-1})
\P(x_i'\gn X_1^{i-1})
\cdot
\\\nonumber&
\ssum{x}{i+1}{n}
[
\P(x_{i+1}^n\gn X_1^{i-1} x_i)\f(X_1^{i-1}x_i' x_{i+1}^n)
-\P(x_{i+1}^n\gn X_1^{i-1} x_i')\f(X_1^{i-1}x_i x_{i+1}^n)
].
\end{align}
}
Let $\pi$ be 
an optimal coupling realizing the infimum in the transportation cost
distance $T_{\rho_{i+1}^n}$
used to define $\tau_{i}(x_1^i,x_i')$.
Recalling (\ref{eq:fg}), we have
\shortlong{
\begin{align}
\nonumber
&
\ssum{x}{i+1}{n}[
\P(x_{i+1}^n\gn X_1^{i-1} x_i')\f(X_1^{i-1}x_i' x_{i+1}^n)
\\\nonumber
-&\P(x_{i+1}^n\gn X_1^{i-1} x_i)\f(X_1^{i-1}x_i x_{i+1}^n)
]
\\\nonumber
=&
~\E
_{(\sseq{\dot X}{i+1}n,\sseq{\ddot X}{i+1}n)\sim\pi}
\sqprn{
\f(\sseq{X}1{i-1}x_i\sseq{\dot X}{i+1}n)
-
\f(\sseq{X}1{i-1}x_i'\sseq{\ddot X}{i+1}n)
}\\\nonumber
\le&
\E
_{
\pi}
\!
\big[
\f(\sseq{X}1{i-1}x_i\sseq{\dot X}{i+1}n)
\!-\!
\f(\sseq{X}1{i-1}x_i'\sseq{\dot X}{i+1}n)
\!+\!
\narsum{j=i+1}^n \rho_j(\dot X_j,\!\ddot X_j)
\big]
\\\nonumber\le&
\E
_{
\sseq{\dot X}{i+1}n
\!\sim \P(\cdot| X_1^{i-1}\!x_i)
}
\!\big[
\f(\sseq{X}1{i-1}x_i\sseq{\dot X}{i+1}n\!)
\!-\!
\f(\sseq{X}1{i-1}x_i'\sseq{\dot X}{i+1}n\!)
\big]
\!\!+\!
\bar\tau_i
\\=&
\label{eq:Ftau}
~F(x_i)-F(x_i')+
\bar\tau_i
,
\end{align}
}{
\begin{align}
\nonumber
&
\ssum{x}{i+1}{n}\sqprn{
\P(x_{i+1}^n\gn X_1^{i-1} x_i)\f(X_1^{i-1}x_i x_{i+1}^n)
-\P(x_{i+1}^n\gn X_1^{i-1} x_i')\f(X_1^{i-1}x_i' x_{i+1}^n)
}
\\\nonumber
=&
~\E
_{(\sseq{\dot X}{i+1}n,\sseq{\ddot X}{i+1}n)\sim\pi}
\sqprn{
\f(\sseq{X}1{i-1}x_i\sseq{\dot X}{i+1}n)
-
\f(\sseq{X}1{i-1}x_i'\sseq{\ddot X}{i+1}n)
}\\\nonumber
\le&
~\E
_{
(\sseq{\dot X}{i+1}n,\sseq{\ddot X}{i+1}n)\sim
\pi}
\sqprn{
\f(\sseq{X}1{i-1}x_i\sseq{\dot X}{i+1}n)
-
\f(\sseq{X}1{i-1}x_i'\sseq{\dot X}{i+1}n)
+
\narsum{j=i+1}^n \rho_j(\dot X_j,\ddot X_j)
}
\\\nonumber\le&
~\E
_{
\sseq{\dot X}{i+1}n
\sim \P(\cdot\gn X_1^{i-1}x_i)
}
\sqprn{
\f(\sseq{X}1{i-1}x_i\sseq{\dot X}{i+1}n)
-
\f(\sseq{X}1{i-1}x_i'\sseq{\dot X}{i+1}n)
}
+
\bar\tau_i
\\=&
\label{eq:Ftau}
~F(x_i)-F(x_i')+
\bar\tau_i
,
\end{align}
}where 
the first inequality holds
by the Lipschitz property and the second by definition of $\bar\tau_{i}$,
and $F:\X_i\to\R$ is defined by
\beq
F(y) = 
\ssum{x}{i+1}{n}\P( \sseq{x}{i+1}{n}\gn X_1^{i-1}x_i) 
\f(\sseq{X}{1}{i-1}y\sseq{x}{i+1}n).
\eeq
Let us substitute (\ref{eq:Ftau}) into (\ref{eq:Vtil}):
\begin{align*}
\tilde V_i&\le
\bar\tau_i
+
\sum_{x_i,x_i'} 
\P( x_i\gn X_1^{i-1})
\P(x_i'\gn X_1^{i-1})
(F(x_i)-F(x_i')).
\end{align*}
Observe that $F$ is $1$-Lipschitz under $\rho_i$ and
apply Jensen's inequality:
\shortlong{
\begin{align*}
&\E[e^{\lam  V_i}\gn X_1^{i-1}]\\
&\le e^{\lam\bar\tau_i}
\sum_{x_i,x_i'} 
\P( x_i\gn X_1^{i-1})
\P(x_i'\gn X_1^{i-1})
e^{\lam(F(x_i)-F(x_i'))}\\
&\le e^{\lam\bar\tau_i}
\sum_{x_i,x_i'} 
\P( x_i\gn X_1^{i-1})
\P(x_i'\gn X_1^{i-1})
\cosh(\lam\rho(x_i,x_i'))
\\
&\le
\exp\paren{\lam
\bar\tau_i
+
\oo2 \bsubg^2(\X_i)\lam^2}
,
\end{align*}
}{
\begin{align*}
\E[e^{\lam  V_i}\gn X_1^{i-1}]
&\le e^{\lam\bar\tau_i}
\sum_{x_i,x_i'} 
\P( x_i\gn X_1^{i-1})
\P(x_i'\gn X_1^{i-1})
e^{\lam(F(x_i)-F(x_i'))}\\
&\le e^{\lam\bar\tau_i}
\sum_{x_i,x_i'} 
\P( x_i\gn X_1^{i-1})
\P(x_i'\gn X_1^{i-1})
\cosh(\lam\rho(x_i,x_i'))
\\
&\le
\exp\paren{\lam
\bar\tau_i
+
\oo2 \bsubg^2(\X_i)\lam^2}
,
\end{align*}
}where the second inequality follows from the argument in (\ref{eq:symm-arg})
and the 
third from
the definition 
of $\bsubg(\X_i)$.
Repeating the standard martingale argument in (\ref{eq:exp-bd}) yields
\shortlong{
\begin{align*}
&\Pr{ \f - \E \f > t} =
\Pr{ \sum_{i=1}^n V_i  > t} 
\\&\le
\exp\paren{\oo2\lam^2\sum_{i=1}^n\bsubg^2(\X_i)-\lam t + 
\lam\sum_{i=1}^n
\bar\tau_i
}.
\end{align*}
}{
\begin{align*}
\Pr{ \f - \E \f > t} &=
\Pr{ \sum_{i=1}^n V_i  > t} 
\\&\le
\exp\paren{\oo2\lam^2\sum_{i=1}^n\bsubg^2(\X_i)-\lam t + 
\lam\sum_{i=1}^n
\bar\tau_i
}.
\end{align*}
}Optimizing over $\lam$ yields the claim.
\enpf

\section{Other Orlicz diameters}
\label{sec:other-orlicz}
Let us recall the notion of an {\em Orlicz norm} $\nrm{X}_\Psi$ of a real random variable $X$ (see, e.g., \citet{MR1113700}):
\beq
\nrm{X}_\Psi = \inf\set{
t>0 : \E[\Psi(X/t)] \le1
},
\eeq
where $\Psi:\R\to\R$ is a {\em Young function} --- i.e., nonnegative, even, convex and vanishing at $0$.
In this section, we will consider the Young functions
\beq
\psi_p(x) = e^{|x|^p}-1,
\qquad p>1,
\eeq
and their induced Orlicz norms. 
A random variable $X$ is subgaussian if and only if $\nrm{X}_{\psi_2}<\infty$.
For $p\neq 2$, $\nrm{X}_{\psi_p}<\infty$ implies that
\beqn
\label{eq:p-orl}
\E e^{\lam X} \le e^{(a|\lam|)^p/p},\qquad \lam\in\R,
\eeqn
for some $a>0$, but the converse implication need not hold.
An immediate consequence of 
Markov's inequality is that any $X$ for which (\ref{eq:p-orl}) holds also satisfies
\beqn
\label{eq:p-orl-conc}
\P(|X|\ge t)\le 2\exp\paren{-\frac{p-1}{p}\paren{\frac{t}{a}}^{p/(p-1)}}.
\eeqn
We define the $p$-{\em Orlicz diameter} of a metric probability space $(\X,\rho,\mu)$,
denoted $\orlp(\X)$,
as the smallest $a>0$ that verifies (\ref{eq:p-orl}) for the symmetrized distance $\Xi(\X)$.
In light of (\ref{eq:p-orl-conc}),
Theorem~\ref{thm:subg-indep} extends straightforwardly to finite $p$-Orlicz metric diameters:

\begin{theorem}
\label{thm:p-orl}
Let $(\X_i,\rho_i,\mu_i)$, $i=1,\ldots,n$  
be a 
sequence of
metric probability spaces and equip $\X^n$ with the usual product measure $\mu^n$ and $\ell_1$ product metric $\rho^n$.
Suppose that for some $p>1$ and all $i\in[n]$ we have $\orlp(\X_i)<\infty$,
and define the vector $\bs{\Delta}\in\R^n$ by $\Delta_i=\orlp(\X_i)$.
If 
$\f:\X\prd\to\R$ is $1$-Lipschitz then for all $t>0$,
\beq
\Pr{ 
\abs{
\f
- 
\E \f
}
> t} \le 
2\exp\paren{
-\frac{p-1}{p}
\paren{\frac{t}{\nrm{\bs{\Delta}}_p}}^{p/(p-1)}
}.
\eeq
\end{theorem}

\section{Discussion}
\label{sec:open-prob}
We have given a concentration inequality for metric spaces with unbounded diameter,
showed its applicability to algorithmic stability with unbounded losses,
and gave an extension to non-independent sampling processes.
Some fascinating questions remain:
\begin{itemize}
\item[(i)] How tight is Theorem~\ref{thm:subg-indep}?
First there is the vexing matter of having a worse constant
in the exponent (i.e., $1/2$) than McDiarmid's (optimal) constant $2$.
Although this gap is not of critical importance, one would like a bound
that recovers McDiarmid's in the finite-diameter case.
More importantly,
is it the case that finite subgaussian diameter is necessary
for subgaussian concentration of all Lipschitz functions?
That is,
given the metric probability spaces $(\X_i,\rho_i,\mu_i)$, $i\in[n]$,
can one always exhibit a $1$-Lipschitz $\f:\X^n\to\R$ that
achieves a nearly matching lower bound?
\item[(ii)] We would like to better understand how 
Theorem~\ref{thm:subg-indep} compares to the Kutin-Niyogi bound (\ref{eq:k-n-ineq}).
We conjecture that for any $(\X^n,\mu^n)$ and $\f:\X^n\to\R$ that satisfies
(\ref{eq:k-n-b}) and (\ref{eq:k-n-c}), one can construct a product metric $\rho^n$
for which $\sum_{i\in[n]}\subg^2(\X_i)<\infty$ and $\f$ is $1$-Lipschitz.
This would imply that whenever the Kutin-Niyogi bound is nontrivial,
so is Theorem~\ref{thm:subg-indep}. We have already shown by example (\ref{eq:ex1})
that the reverse does not hold.
\item[(iii)] 
The quantity $\bar\tau_i$ defined in Section~\ref{sec:mix} is a rather complicated object;
one desires a better handle on it in terms of the given distribution and metric.
\item[(iv)]
Perhaps the most pressing question is that of showing that some common
learning algorithms such as $k$-nearest neighbor, kernel SVM, regularized regression
are totally Lipschitz stable under our definition (\ref{eq:tot-lip}).
\end{itemize}

\subsubsection*{Acknowledgements}
\shortlong{
We thank some folks.
}{
John Lafferty encouraged me to seek a distribution-dependent refinement of McDiarmid's inequality.
Thanks also to Gideon Schechtman, Shahar Mendelson, Assaf Naor, Iosif Pinelis 
and Csaba Szepesv\'ari
for helpful correspondence, and to Roi Weiss for carefully proofreading the manuscript.
}

\small
\bibliography{../../mybib}
\bibliographystyle{plainnat}

\end{document}